\documentclass[12pt]{article}
\usepackage{txfonts}

\def\T{\Gamma} \def\D{\Delta} \def\Th{\Theta}
\def\Ld{\Lambda} \def\E{\Sigma} \def\O{\Omega}
\def\a{\alpha} \def\b{\beta} \def\g{\gamma} \def\d{\delta} \def\e{\varepsilon}
\def\r{\rho} \def\o{\sigma} \def\t{\tau} \def\w{\omega} \def\k{\kappa}
\def\th{\theta} \def\ld{\lambda} \def\ph{\varphi} \def\z{\zeta}
\def\A{$A$~} \def\G{$G$~} \def\H{$H$~} \def\K{$K$~} \def\M{$M$~} \def\N{$N$~}
\def\P{$P$~} \def\Q{$Q$~} \def\R{$R$~} \def\V{$V$~} \def\X{$X$~} \def\Y{$Y$~}
\def\rmA{{\bf A}} \def\rmD{{\bf D}} \def\rmS{{\bf S}} \def\rmK{{\bf K}}
\def\rmM{{\bf M}} \def\Z{{\Bbb Z}} \def\GL{{\bf GL}} \def\C{Cayley}
\def\oa{\ovl A} \def\og{\ovl G} \def\oh{\ovl H} \def\ob{\ovl B} \def\oq{\ovl Q}
\def\oc{\ovl C} \def\ok{\ovl K} \def\ol{\ovl L} \def\om{\ovl M} \def\on{\ovl N}
\def\op{\ovl P} \def\oR{\ovl R} \def\os{\ovl S} \def\ot{\ovl T} \def\ou{\ovl U}
\def\ov{\ovl V} \def\ow{\ovl W} \def\ox{\ovl X} \def\oT{\ovl\T}
\def\lg{\langle} \def\rg{\rangle}
\def\di{\bigm|} \def\Di{\Bigm|} \def\nd{\mathrel{\bigm|\kern-.7em/}}
\def\Nd{\mathrel{\not\,\Bigm|}} \def\edi{\bigm|\bigm|}
\def\m{\medskip} \def\l{\noindent} \def\x{$\!\,$}  \def\J{$-\!\,$}
\def\Hom{\hbox{\rm Hom}} \def\Aut{\hbox{\rm Aut}} \def\Inn{\hbox{\rm Inn}}
\def\Syl{\hbox{\rm Syl}} \def\Sym{\hbox{\rm Sym}} \def\Alt{\hbox{\rm Alt}}
\def\Ker{\hbox{\rm Ker}} \def\fix{\hbox{\rm fix}} \def\mod{\hbox{\rm mod}}
\def\psl{{\bf P\!SL}} \def\Cay{\hbox{\rm Cay}} \def\Mult{\hbox{\rm Mult}}
\def\val{\hbox{\rm Val}} \def\Sab{\hbox{\rm Sab}} \def\supp{\hbox{\rm supp}}
\def\qed{\hfill $\Box$} \def\qqed{\qed\vspace{3truemm}}
\def\CS{\Cay(G,S)} \def\CT{\Cay(G,T)}
\def\h{\heiti\bf} \def\hs{\ziti{E}\bf} \def\st{\songti} \def\ft{\fangsong}
\def\kt{\kaishu} \def\heit{\hs\relax} \def\songt{\st\rm\relax}
\def\fangs{\ft\rm\relax} \def\kaish{\kt\rm\relax} \def\fs{\footnotesize}
\begin{document}

\newtheorem{theorem}{Theorem}[section]
\newtheorem{corollary}[theorem]{Corollary}
\newtheorem{definition}[theorem]{Definition}
\newtheorem{conjecture}[theorem]{Conjecture}
\newtheorem{question}[theorem]{Question}
\newtheorem{lemma}[theorem]{Lemma}
\newtheorem{proposition}[theorem]{Proposition}
\newtheorem{example}[theorem]{Example}
\newtheorem{problem}[theorem]{Problem}
\newenvironment{proof}{\noindent {\bf
Proof.}}{\rule{3mm}{3mm}\par\medskip}
\newcommand{\remark}{\medskip\par\noindent {\bf Remark.~~}}
\newcommand{\pp}{{\it p.}}
\newcommand{\de}{\em}

%=================================================
\def\dfrac{\displaystyle\frac} \def\ovl{\overline}
\def\for{\forall~} \def\exi{\exists~} \def\c{\subseteq}
\def\iif{\Longleftrightarrow} \def\Rto{\Rightarrow} \def\Lto{\Leftarrow}
%=================================================
\def\T{\Gamma} \def\D{\Delta} \def\Th{\Theta}
\def\Ld{\Lambda} \def\E{\Sigma} \def\O{\Omega}
\def\a{\alpha} \def\b{\beta} \def\g{\gamma} \def\d{\delta} \def\e{\varepsilon}
\def\r{\rho} \def\o{\sigma} \def\t{\tau} \def\w{\omega} \def\k{\kappa}
\def\th{\theta} \def\ld{\lambda} \def\ph{\varphi} \def\z{\zeta}
%=================================================
\def\A{$A$~} \def\G{$G$~} \def\H{$H$~} \def\K{$K$~} \def\M{$M$~} \def\N{$N$~}
\def\P{$P$~} \def\Q{$Q$~} \def\R{$R$~} \def\V{$V$~} \def\X{$X$~} \def\Y{$Y$~}
\def\rmA{{\bf A}} \def\rmD{{\bf D}} \def\rmS{{\bf S}} \def\rmK{{\bf K}}
\def\rmM{{\bf M}} \def\Z{{\Bbb Z}} \def\GL{{\bf GL}} \def\C{Cayley}
%=================================================
\def\oa{\ovl A} \def\og{\ovl G} \def\oh{\ovl H} \def\ob{\ovl B} \def\oq{\ovl Q}
\def\oc{\ovl C} \def\ok{\ovl K} \def\ol{\ovl L} \def\om{\ovl M} \def\on{\ovl N}
\def\op{\ovl P} \def\oR{\ovl R} \def\os{\ovl S} \def\ot{\ovl T} \def\ou{\ovl U}
\def\ov{\ovl V} \def\ow{\ovl W} \def\ox{\ovl X} \def\oT{\ovl\T}
\def\lg{\langle} \def\rg{\rangle}
%=================================================
\def\di{\bigm|} \def\Di{\Bigm|} \def\nd{\mathrel{\bigm|\kern-.7em/}}
\def\Nd{\mathrel{\not\,\Bigm|}} \def\edi{\bigm|\bigm|}
\def\m{\medskip} \def\l{\noindent} \def\x{$\!\,$}  \def\J{$-\!\,$}
%=================================================
\def\Hom{\hbox{\rm Hom}} \def\Aut{\hbox{\rm Aut}} \def\Inn{\hbox{\rm Inn}}
\def\Syl{\hbox{\rm Syl}} \def\Sym{\hbox{\rm Sym}} \def\Alt{\hbox{\rm Alt}}
\def\Ker{\hbox{\rm Ker}} \def\fix{\hbox{\rm fix}} \def\mod{\hbox{\rm mod}}
\def\psl{{\bf P\!SL}} \def\Cay{\hbox{\rm Cay}} \def\Mult{\hbox{\rm Mult}}
\def\val{\hbox{\rm Val}} \def\Sab{\hbox{\rm Sab}} \def\supp{\hbox{\rm supp}}
\def\qed{\hfill $\Box$} \def\qqed{\qed\vspace{3truemm}}
\def\CS{\Cay(G,S)} \def\CT{\Cay(G,T)}
%=================================================
\def\h{\heiti\bf} \def\hs{\ziti{E}\bf} \def\st{\songti} \def\ft{\fangsong}
\def\kt{\kaishu} \def\heit{\hs\relax} \def\songt{\st\rm\relax}
\def\fangs{\ft\rm\relax} \def\kaish{\kt\rm\relax} \def\fs{\footnotesize}
%=================================================

\newcommand{\JEC}{{\it Europ. J. Combinatorics},  }
\newcommand{\JCTB}{{\it J. Combin. Theory Ser. B.}, }
\newcommand{\JCT}{{\it J. Combin. Theory}, }
\newcommand{\JGT}{{\it J. Graph Theory}, }
\newcommand{\ComHung}{{\it Combinatorica}, }
\newcommand{\DM}{{\it Discrete Math.}, }
\newcommand{\ARS}{{\it Ars Combin.}, }
\newcommand{\SIAMDM}{{\it SIAM J. Discrete Math.}, }
\newcommand{\SIAMADM}{{\it SIAM J. Algebraic Discrete Methods}, }
\newcommand{\SIAMC}{{\it SIAM J. Comput.}, }
\newcommand{\ConAMS}{{\it Contemp. Math. AMS}, }
\newcommand{\TransAMS}{{\it Trans. Amer. Math. Soc.}, }
\newcommand{\AnDM}{{\it Ann. Discrete Math.}, }
\newcommand{\NBS}{{\it J. Res. Nat. Bur. Standards} {\rm B}, }
\newcommand{\ConNum}{{\it Congr. Numer.}, }
\newcommand{\CJM}{{\it Canad. J. Math.}, }
\newcommand{\JLMS}{{\it J. London Math. Soc.}, }
\newcommand{\PLMS}{{\it Proc. London Math. Soc.}, }
\newcommand{\PAMS}{{\it Proc. Amer. Math. Soc.}, }
\newcommand{\JCMCC}{{\it J. Combin. Math. Combin. Comput.}, }
\newcommand{\GC}{{\it Graphs Combin.}, }

\title{ Automorphism Groups of the Pancake Graphs \thanks{
This work is supported by National Natural Science Foundation of
China (No:10971137), the National Basic Research Program (973) of
China (No.2006CB805900),  and a grant of Science and Technology
Commission of Shanghai Municipality (STCSM, No: 09XD1402500)
.\newline \indent
 $^{\dagger}$Correspondent author: Xiao-Dong
Zhang (Email: xiaodong@sjtu.edu.cn)}}
\author{
 Yun-Ping Deng, Xiao-Dong Zhang$^{\dagger}$\\
{\small Department of Mathematics,
 Shanghai Jiao Tong University} \\
{\small  800 Dongchuan road, Shanghai, 200240, P.R. China}\\
{\small Emails: dyp612@hotmail.com, xiaodong@sjtu.edu.cn}\\
 }
\date{}
\maketitle
 \begin{abstract}
   It is well-known that the pancake graphs are widely used as models for interconnection networks \cite{Akers}.
   In this paper, some properties of the pancake graphs are investigated. We first prove that the pancake graph, denoted by $P_n~(n\geq 4),$ is
   super-connected and hyper-connected. Further,
   we study the symmetry of $P_n$ and completely determine its full automorphism group,
    which shows that $P_n~(n\geq 5)$ is
   a graphical regular representation of $S_n.$

 \end{abstract}

{{\bf Key words:} Interconnection networks; pancake graph; super-connected; hyper-connected; efficient dominating sets; automorphism group. }\\

{{\bf AMS Classifications:} 05C25, 05C69}
\vskip 0.5cm

\section{Introduction}
For a simple graph $\T,$ we denote its vertex set, edge set and full
automorphism group respectively by $V(\T), E(\T)$ and $\Aut(\T)$.
$\T$ is said to be {\it vertex-transitive} or {\it edge-transitive}
if $\Aut(\T)$ acts transitively on $V(\T)$ or $E(\T),$ respectively.
Let \G be a finite group and $S$ a subset of \G not containing the
identity element $1$ with $S=S^{-1}.$ The {\it Cayley graph}
$\T:=\Cay(G,S)$ on \G with respect to $S$ is defined by
$$V(\T){=}G,~E(\T){=}\{(g,gs): g{\in} G,\ s{\in} S\}.$$ Clearly, $\T$
is a $|S|$-regular and vertex-transitive graph, since $\Aut(\T)$
contains the left regular representation $L(G)$ of \G. Moreover,
$\T$ is connected if and only if \G is generated by $S.$

A {\it permutation} $\sigma$ on the set $X=\{1,2,\cdots,n\}$ is a
bijective mapping from $X$ to $X.$ As usual, we denote by $S_n$ the
group of all permutations on $X,$ which is called the {\it symmetric
group}. The {\it pancake graph} $P_n,$ also called the {\it
prefix-reversal graph} is the Cayley graph $\Cay(S_n,PR_n),$ where
$PR_n= \{{r_{1j}:2\leq j\leq n}\}$ and $r_{1j}=(\begin{array}{ll}
1~~~~2~~\,\cdots~j~j+1~\cdots~n\\
j~j-1\cdots~1~j+1~\cdots~n
\end{array}).$

The pancake graph is well-known because of the famous unsolved
combinatorial problem about computing its diameter, which has been
introduced by \cite{Dweighter}, and has been studied in several
papers\cite{Gates,Hyedari1,Hyedari2}. The pancake graph was often
used as a model for interconnection networks of parallel computers
\cite{Akers} due to its attractive properties regarding degree,
diameter, symmetry, embeddings and self similarity. The pancake
graph $P_n$ corresponds to the $n$-dimensional pancake network in
computer science such that this network has processors labeled by
permutations on $X$ and two processors are connected when the label
of one is obtained from the other by some $r_{1j}, 2\leq j\leq n.$
The diameter of this network corresponds to the worst communication
delay for transmitting information in a system. Morover, many
researchers (see \cite{Kanevsky,Lin,Sheu}) have investigated some
other properties of $P_n,$ such as the hamilton-connectedness,
cycle-embedding problem, super-connectivity.

A graph \X is said to be {\it super-connected} \cite{Boesch} if each
minimum vertex cut is the neighbor set of a single vertex in \X. A
graph \X is said to be {\it hyper-connected} \cite{Hamidoune} if for
every minimum vertex cut $D$ of \X, $X-D$ has exactly two
components, one of which is an isolated vertex. In \cite{Li}, Li
investigated the super-connectedness and hyper-connectedness of the
reversal Cayley graph and pointed out that it is unknown for the
pancake graph. Here we solve this problem and prove that the pancake
graph $P_n~(n\geq 4)$ is super-connected and hyper-connected.

An independent set $D$ of vertices in a graph is called an {\it
efficient dominating set} \cite{Bange1,Bange2} if each vertex not in
$D$ is adjacent to exactly one vertex in $D.$ In \cite{Dejter},
Dejter investigated the efficient dominating sets of Cayley graphs
on the symmetric groups, which implied that there exists the
efficient dominating sets in the pancake graph. In addtion, the
efficient dominating sets are used in optimal broadcasting
algorithms for multiple messages on the pancake graphs \cite{Qiu}.
Motivated by these results, we completely characterize all the
efficient dominating sets in $P_n~(n\geq 3).$

A graph $\T=(V,E)$ is a {\it graphical regular representation
(GRR)}\cite{Nowitz} of the finite group $G$ if $\Aut(\T)=G$ and
$\Aut(\T)$ acts regularly on $V.$ It is well-known that for the
interconnection networks modeled by Cayley graphs, the symmetry is
one of the problems focused by many researchers. In
\cite{Lakshmivarahan}, Lakshmivarahan investigated the symmetry of
the pancake graph and showed that $P_n$ is not edge-transitive and
hence not distance-transitive. In this paper, we further study the
symmetry of $P_n$ and completely determine the automorphism group of
$P_n,$ which shows that $P_n~(n\geq 5)$ is a graphical regular
representation of $S_n$ and hence not edge-transitive and
distance-transitive.

The rest part of this paper is organized as follows. In Section 2,
we first prove that the pancake graph $P_n~(n\geq 4)$ is
super-connected and hyper-connected, then we show that there are
exactly $n$ efficient dominating sets $B^{(i)}~(i=1,2,\cdots,n)$ in
$P_n~(n\geq 3),$ where $B^{(i)}=\{\pi\in S_n:\,\pi(1)=i\}.$ In
section 3, we prove that the full automorphism groups of
$P_n~(n\geq5)$ is the left regular representation of $S_n,$ i.e.
$\Aut(P_n)=L(S_n).$

\section{Some properties of $P_n$}

 In table 1 of \cite {Li}, it has been pointed out that the super-connectedness
 and hyper-connectedness of $P_n$ are unknown.  In this section, we first prove
 that the pancake graph $P_n~(n\geq 4)$ is  super-connected and hyper-connected.
 Following \cite {Li}, we introduce  some notations and terminologies. Let $X$
 be a graph and $F$ a subset of  $V(X).$ Set $N(F)=\{x\in V(X)\setminus F:\exists\,y\in
 F,\,s.t.\,xy\in E(X)\},\,C(F)=F\cup N(F),\,R(F)=V(X)\setminus
 C(F).$ A subset $F\subseteq V(X)$ is a {\it fragment} if  $|N(F)|=\kappa(X)$ and
 $R(F)\neq\emptyset,$ where $\kappa(X)$ is the vertex-connectivity of $X.$ A fragment $F$ with
 $2\leq|F|\leq|V(X)|-\kappa(X)-2$ is called a {\it strict fragment}.  A strict fragment with
 minimum cardinality is called a {\it superatom}.

The following result is due to Mader \cite{Mader}:

\begin{lemma} ${{\fs\cite{Mader}}}$ \label{yl-2.1}
If $X$ is a connected undirected graph which is a vertex-transitive
and $K_4$-free, then $\kappa(X)=\delta(X),$ where $\delta(X)$
denotes the minimum degree of $X.$
\end{lemma}

\begin{lemma} \label{yl-2.2}
$\kappa(P_n)=\delta(P_n)=n-1$ for any $n\geq 3.$
\end{lemma}

\begin{proof} By \cite{Sheu}, we obtain that $g(P_n)=6,$ where $g(P_n)$ is the girth of $P_n.$ So
$P_n$ is $K_4$-free, by Lemma \ref{yl-2.1}, the assertion
holds.\end{proof}

In the following Lemma, we shall state some facts without proof.
Some of these facts may be found in \cite {Sheu}, and others follow
immediately from the definition of the pancake graph.

\begin{lemma} \label{yl-2.3}
Let $B^{(i)}=\{\pi\in S_n:\,\pi(1)=i\},\,B_{(j)}=\{\pi\in
S_n:\,\pi(n)=j\},\,B_{(j)}^{(i)}=B^{(i)}\cap B_{(j)}.$ Then the
following (i)-(iii) hold:

(i) For any $i\neq j,$ each vertex in $B^{(i)}$ is adjacent to
exactly one vertex in $B^{(j)};$

(ii) For any $i\neq j,$ each vertex in $B_{(j)}^{(i)}$ is adjacent
to exactly one vertex in $B_{(i)}^{(j)}$ and exactly one vertex in
$B_{(j)}^{(k)}$ for each $k\neq i,j.$

(iii) The mapping $\varphi:\,S_{n-1}\rightarrow B_{(j)}$ defined as
$\varphi(\pi)=(j,n)\pi$ is an isomorphism from $P_{n-1}$ to
$P_n[B_{(j)}],$ where $P_n[B_{(j)}]$ is the subgraph of $P_n$
induced by $B_{(j)}.$
\end{lemma}

\begin{theorem} \label{dl-2.4}
If $n\geq 4,$ then $P_n$ is super-connected.
\end{theorem}

\begin{proof} It is enough to show that $P_n$ contains no superatom.
Suppose on the contrary that $A$ is a superatom of $P_n$ and
consider the following possible cases:

{\bf Case 1.} $A\subseteq B_{(i)}$ for some $i\in\{1,2,\cdots,n\}.$

By Lemmas \ref{yl-2.2} and \ref{yl-2.3}, we have
$\kappa(P_{n}[B_{(i)}])=\kappa(P_{n-1})=n-2$ for $n-1\geq 3,$ so
$|N(A)\cap B_{(i)}|\geq n-2$ for $n\geq 4.$ Hence $|N(A)|=|N(A)\cap
B_{(i)}|+|N(A)\cap(\bigcup_{j\neq i} B_{(j)})|\geq
(n-2)+|A|\geq(n-2)+2=n>n-1=\kappa(P_n),$ which is a contradiction.

{\bf Case 2.} $A\nsubseteq B_{(i)}$ for any $i\in\{1,2,\cdots,n\}.$

Then there exist $i,j~(i\neq j)$ such that $A\cap B_{(i)}\neq
\emptyset$ and $A\cap B_{(j)}\neq \emptyset.$ Hence
$|N(A)|\geq|N(A)\cap B_{(i)}|+|N(A)\cap B_{(j)}|\geq 2(n-2)\geq
n-1=\kappa(P_n),$ which is a contradiction.
\end{proof}

{\bf Remark.} If $n=3,$ then $P_3=C_6,$ clearly it is not
super-connected.

\begin{theorem} \label{dl-2.5}
If $n\geq 4,$ then $P_n$ is hyper-connected.
\end{theorem}

\begin{proof} By Theorem \ref{dl-2.4}, $P_n$ is super-connected for $n\geq 4.$ Consider the
vertex-transitivity of $P_n,$ it suffices to show that $P_n-N[I]$ is
connected, where $N[I]$ is the closed neighbourhood of the identity
element $I.$ We proceed by the induction on $n.$ If $n=4,$ one can
easily check that $P_4-N[I]$ is connected.

If $n>4,$ then $P_{n}[B_{(i)}]-N[I]$ is connected for any $i<n$
since $|N[I]\cap B_{(1)}|=1$ and $|N[I]\cap B_{(i)}|=0$ for any
$1<i<n.$ By induction, $P_{n}[B_{(n)}]-N[I]=P_{n-1}-N[I]$ is
connected. By Lemma \ref{yl-2.3}, each vertex in $B_{(i)}^{(n)}$ is
adjacent to exactly one vertex in $B_{(n)}^{(i)}$ for any $i<n.$ So
for each $i<n$ there exists a vertex in $P_{n}[B_{(i)}]-N[I]$ which
is adjacent to some vertex in $P_{n}[B_{(n)}]-N[I].$ Thus
$P_n-N[I]=\bigcup_{i=1}^{n}P_{n}[B_{(i)}]-N[I]$ is connected.
\end{proof}

Next we turn to consider the efficient dominating sets of $P_n.$ By
the definition of efficient dominating set, it is easy to see that
any efficient dominating set $D$ in $P_n$ has $(n-1)!$ elements and
$d(u,v)\geq 3$ for any $u,v\in D,$ where $d(u,v)$ is the distance
between two vertice $u$ and $v$ in $P_n.$ Konstantinova in the
abstract \cite{Elena} obtained the following result on the efficient
dominating set. For the completeness of this paper, here we present
a proof of the result.

\begin{theorem} ${{\fs\cite{Elena}}}$ \label{dl-2.6}
There are exactly $n$ efficient dominating sets $B^{(i)}~(1\leq
i\leq n)$ in $P_n~(n\geq 3).$
\end{theorem}

\begin{proof}
Clearly each $B^{(i)}~(1\leq i\leq n)$ is an efficient dominating
set in $P_n.$ So it suffices to prove that for any efficient
dominating set $D$ in $P_n,$ if $D\cap B^{(i)}\neq\emptyset,$ then
$D=B^{(i)}.$ Set $D^{(i)}=D\cap B^{(i)},\,D_{(j)}=D\cap
B_{(j)},\,D_{(j)}^{(i)}=D\cap
B^{(i)}_{(j)},\,R_{(j)}^{(i)}=B^{(i)}_{(j)}\setminus D_{(j)}^{(i)}.$
We consider the following cases:

{\bf Case 1.} There exists a $D_{(j)}^{(i)}$ such that
$D_{(j)}^{(i)}=B^{(i)}_{(j)}.$

By Lemma \ref{yl-2.3}, $N(B^{(i)}_{(j)})\cap
B_{(i)}=B^{(j)}_{(i)},\,N(N(B^{(i)}_{(j)}))\cap
B_{(i)}=B_{(i)}\setminus B^{(j)}_{(i)}.$ Since
$B^{(i)}_{(j)}\subseteq D$ and $d(u,v)\geq 3$ for any $u,v\in D,$ so
we have $D\cap B^{(j)}_{(i)} =\emptyset,\,D\cap (B_{(i)}\setminus
B^{(j)}_{(i)})=\emptyset,$ i.e. $D\cap B_{(i)}=\emptyset.$ By Lemma
\ref{yl-2.3} again, $N(B_{(i)})=B^{(i)}$ and each vertex in
$B_{(i)}$ is adjacent to exactly one vertex in $B^{(i)}.$ Hence
$B^{(i)}\subseteq D.$ Since $|B^{(i)}|=|D|=(n-1)!,$ we have
$D=B^{(i)}.$

{\bf Case 2.} There exists a $D_{(j)}^{(i)}$ such that
$\emptyset\neq D_{(j)}^{(i)}\varsubsetneq B^{(i)}_{(j)}.$

By Lemma \ref{yl-2.3}, we have $X_{(i)}:=N(D_{(j)}^{(i)})\cap
B_{(i)}\subseteq B_{(i)}^{(j)},\,Y_{(i)}:=N(R_{(j)}^{(i)})\cap
B_{(i)}\cap D\subseteq
B_{(i)}^{(j)},\,Z_{(i)}:=B_{(i)}^{(j)}\setminus(X_{(i)}\cup
Y_{(i)}),\,W_{(i)}:=N(Z_{(i)})\cap D\subseteq B_{(i)}\setminus
B_{(i)}^{(j)},\,Y_{(j)}:=N(R_{(j)}^{(i)})\cap B_{(j)}\cap D\subseteq
B_{(j)}\setminus B_{(j)}^{(i)}.$ Now we claim that
$D_{(i)}=Y_{(i)}\cup W_{(i)},\,D_{(j)}=D_{(j)}^{(i)}\cup Y_{(j)}.$
Clearly $D_{(i)}\supseteq Y_{(i)}\cup W_{(i)},\,D_{(j)}\supseteq
D_{(j)}^{(i)}\cup Y_{(j)}.$ For any $x\in D_{(i)}\setminus Y_{(i)},$
then $x\in B_{(i)}\setminus (X_{(i)}\cup N(X_{(i)})\cup Y_{(i)}\cup
N(Y_{(i)})\cup Z_{(i)})=N(Z_{(i)})\cap B_{(i)}$ and so $x\in
N(Z_{(i)})\cap B_{(i)}\cap D=W_{(i)}.$ Hence $D_{(i)}\subseteq
Y_{(i)}\cup W_{(i)}.$ For any $y\in D_{(j)}\setminus D_{(j)}^{(i)},$
then $y\in B_{(j)}\setminus (D_{(j)}^{(i)}\cup N(D_{(j)}^{(i)})\cup
R_{(j)}^{(i)})=N(R_{(j)}^{(i)})\cap B_{(j)}$ and so $y\in
N(R_{(j)}^{(i)})\cap B_{(j)}\cap D=Y_{(j)}.$ Hence $D_{(j)}\subseteq
D_{(j)}^{(i)}\cup Y_{(j)}.$

Clearly $|X_{(i)}|=|D_{(j)}^{(i)}|$ and
$|Y_{(i)}|+|Y_{(j)}|=|R_{(j)}^{(i)}|,$ so
$|X_{(i)}|+|Y_{(i)}|+|Y_{(j)}|=|B_{(j)}^{(i)}|=(n-2)!.$ Since
$|X_{(i)}|+|Y_{(i)}|+|Z_{(i)}|=|B_{(i)}^{(j)}|=(n-2)!,$ we have
$|W_{(i)}|=|Z_{(i)}|=|Y_{(j)}|.$ By the definition of efficient
dominating set and Lemma \ref{yl-2.3}, for $k=i,j,$ each vertex in
$B_{(k)}\setminus (D_{(k)}\cup N(D_{(k)})$ is adjacent to exactly
one vertex in $D^{(k)},$ each vertex in $D^{(k)}$ is adjacent to
exactly one vertex in $B_{(k)}\setminus (D_{(k)}\cup N(D_{(k)}).$ So
$|D^{(i)}|=|B_{(i)}\setminus (D_{(i)}\cup N(D_{(i)})|
=(n-1)!-(n-1)|D_{(i)}| =(n-1)!-(n-1)(|Y_{(i)}|+|W_{(i)}|)
=(n-1)!-(n-1)((n-2)!-|X_{(i)}|)=(n-1)|X_{(i)}|,\,
|D^{(j)}|=|B_{(j)}\setminus (D_{(j)}\cup
N(D_{(j)})|=(n-1)!-(n-1)|D_{(j)}|=(n-1)!-(n-1)(|D_{(j)}^{(i)}|+
|Y_{(j)}|)=(n-1)!-(n-1)((n-2)!-|Y_{(i)}|)=(n-1)|Y_{(i)}|.$ Hence
$|\bigcup_{k\neq i,j}D_{(k)}^{(i)}|=|D^{(i)}|-
|D_{(j)}^{(i)}|=(n-1)|X_{(i)}|-|X_{(i)}|=(n-2)|X_{(i)}|,\,|\bigcup_{k\neq
i,j}D_{(k)}^{(j)}|=|D^{(j)}|-|D_{(i)}^{(j)}|=(n-1)|Y_{(i)}|-|Y_{(i)}|=(n-2)|Y_{(i)}|$
and $|\bigcup_{k,l\neq
i,j}D_{(k)}^{(l)}|=|D|-|D^{(i)}|-|D^{(j)}|-|W_{(i)}|-|Y_{(j)}|=
(n-1)!-(n-1)|X_{(i)}|-(n-1)|Y_{(i)}|-|W_{(i)}|-|Y_{(j)}|=
(n-1)((n-2)!-|X_{(i)}|-|Y_{(i)}|)-2|Z_{(i)}|=(n-3)|Z_{(i)}|.$

By the definition of efficient dominating set and Lemma \ref{yl-2.3}
, for any a fixed $l_0\neq i,j,$ each vertex in $\bigcup_{k\neq
i,j}B_{(k)}^{(l_0)}$ either belongs to $\bigcup_{k\neq i,j}D_{(k)}$
or is adjacent to exactly one vertex in $\bigcup_{k\neq
i,j}D_{(k)},$ each vertex in $\bigcup_{k\neq i,j}D_{(k)}$ either
belongs to $\bigcup_{k\neq i,j}B_{(k)}^{(l_0)}$ or is adjacent to
exactly one vertex in $\bigcup_{k\neq i,j}B_{(k)}^{(l_0)},$so
$(n-3)(n-2)!=|\bigcup_{k\neq i,j}B_{(k)}^{(l_0)}|=|\bigcup_{k\neq
i,j}D_{(k)}|=|\bigcup_{k\neq i,j}D_{(k)}^{(i)}|+|\bigcup_{k\neq
i,j}D_{(k)}^{(j)}|+|\bigcup_{k,l\neq
i,j}D_{(k)}^{(l)}|=(n-2)|X_{(i)}|+(n-2)|Y_{(i)}|+(n-3)|Z_{(i)}|=
(n-3)(|X_{(i)}|+|Y_{(i)}|+|Z_{(i)}|)+|X_{(i)}|+|Y_{(i)}|=(n-3)(n-2)!+|X_{(i)}|+|Y_{(i)}|,$
hence $|D_{(j)}^{(i)}|=|X_{(i)}|=0,$ which is a contradiction.
\end{proof}

\section{The automorphism group of $P_n$}

In this section, we completely determine the full automorphism group
of $P_n.$ First we introduce some definitions. Let $\Sym(\Omega)$
denote the set of all permutations of a set $\Omega.$ A {\it
permutation representation} of a group \G is a homomorphism from \G
into $\Sym(\Omega)$ for some set $\Omega.$ A permutation
representation is also referred to as an action of \G on the set
$\Omega,$ in which case we say that \G acts on $\Omega.$
Furthermore, if $\{g\in G:x^g=x,\,\forall x\in \Omega\}=1,$ we say
the action of $G$ on $\Omega$ is {\it faithful}, or $G$ acts {\it
faithfully} on $\Omega.$

\begin{theorem}\label{dl-3.1}
For $n\geq 5,$ if $N(X)=B^{(i)}$ and $|X|=|B^{(i)}|,$ where
$X\subseteq V(P_n)=S_n$ and $i\in \{1,2,\cdots,n\},$ then
$X=B_{(i)}.$
\end{theorem}

\begin{proof}
For $n=5,$ one can easily check that the assertion holds. We proceed
by induction on $n.$ First since $N(X)=\{y\in V(P_n)\setminus
X:\exists\,x\in X,\,s.t.\,xy\in E(P_n)\},$ we have $X\cap
N(X)=\emptyset,$ i.e. $X\cap B^{(i)}=\emptyset.$ Next we shall show
that $X=B_{(i)}$ by the following three Claims:

{\bf Claim 1.} Either $X=B_{(i)}$ or $X\cap B_{(i)}=\emptyset.$

Set $X_i:=X\cap B_{(i)},\,\overline{X_i}:=X\setminus X_i.$ Suppose
on the contrary that $\emptyset\neq X_i\subsetneq B_{(i)}.$ By Lemma
\ref{yl-2.3} (iii), $P_n[B_{(i)}]\cong P_{n-1},$ so $P_n[B_{(i)}]$
is connected, which implies that $N(X_i)\cap B_{(i)}\neq\emptyset.$
Since $N(X_i)\subseteq B_{(i)}\cup B^{(i)}$ and $\overline{X_i}\cap
(B_{(i)}\cup B^{(i)})=\emptyset,$ we have $N(X_i)\cap
\overline{X_i}=\emptyset,$ i.e. $N(X_i)\subseteq N(X).$ So $N(X)\cap
B_{(i)}\neq\emptyset,$ which contradicts $N(X)=B^{(i)},$ hence Claim
1 holds.

{\bf Claim 2.} Set $X_k=X\cap B_{(k)},\,B_{(n-1)\rightarrow
i,\,n\rightarrow k}=\{\pi\in S_n:\,\pi(n-1)=i,\,\pi(n)=k\}.$ If
$X\neq B_{(i)},$ then $X_k=B_{(n-1)\rightarrow i,\,n\rightarrow k}$
for any $k\neq i.$

By $X\neq B_{(i)}$ and Claim 1, $X\cap (B^{(i)}\cup
B_{(i)})=\emptyset.$ By Lemma \ref{yl-2.3} (ii), $B^{(i)}_{(k)}\cap
N(X_l)=\emptyset$ for any $k\neq l.$ So we have
$B^{(i)}_{(k)}\subseteq B^{(i)}=N(X)=N(\bigcup_{k\neq
i}X_k)\subseteq\bigcup_{k\neq i}N(X_k)\Rightarrow
B^{(i)}_{(k)}\subseteq N(X_k)\Rightarrow B^{(i)}_{(k)}\subseteq
B_{(k)}\cap N(X_k).$ On the other hand, $B_{(k)}\cap N(X_k)\subseteq
B_{(k)}\cap N(X)=B_{(k)}\cap B^{(i)}=B^{(i)}_{(k)}.$ Thus
$B_{(k)}\cap N(X_k)=B^{(i)}_{(k)}.$ By Theorem \ref{dl-2.6},
$B^{(i)}$ is an efficient dominating set of $P_n,$ so $|X_k|\geq
|B^{(i)}_{(k)}|\Rightarrow|X|=\sum_{k\neq i}|X_k|\geq\sum_{k\neq
i}|B^{(i)}_{(k)}|=|B^{(i)}|,$ note that $|X|=|B^{(i)}|,$ and so
$|X_k|=|B^{(i)}_{(k)}|.$ By Lemma \ref{yl-2.3} (iii),
$P_n[B_{(k)}]\cong P_{n-1}$ and $B^{(i)}_{(k)}$ is an efficient
dominating set of $P_n[B_{(k)}].$ Since $B_{(k)}\cap
N(X_k)=B^{(i)}_{(k)}$ and $|X_k|=|B^{(i)}_{(k)}|,$ by induction, we
have $X_k=B_{(n-1)\rightarrow i,\,n\rightarrow k},$ hence Claim 2
holds.

{\bf Claim 3.} If $X\neq B_{(i)},$ then $n=3.$

By $X\neq B_{(i)}$ and Claim 2, $X_k=B_{(n-1)\rightarrow
i,\,n\rightarrow k}$ for any $k\neq i.$ Since $X_k\subseteq B_{(k)}$
and $P_n[B_{(k)}]\cong P_{n-1},$ which is a $(n-2)$-regular graph,
we have $|N(x_k)\cap B_{(k)}|=n-2$ for any $x_k\in X_k,$ note that
$|N(x_k)|=n-1,$ and so $|N(x_k)\cap(\bigcup_{l\neq k} B_{(l)})|=1.$
Set $N(x_k)\cap(\bigcup_{l\neq k} B_{(l)})=\{x_l\},$ where $x_l\in
B_{(l)}$ for some $l\neq k,i.$ Since $x_l\in N(x_k)\cap B_{(l)}$ and
$N(x_k)\cap B_{(l)}\cap B^{(i)}=\emptyset$ (by Lemma \ref{yl-2.3}),
we have $x_l\not\in B^{(i)}.$ Note that $x_l\in N(x_k)$ and
$N(X)=B^{(i)},$ then $x_l\in X_l=B_{(n-1)\rightarrow
i,\,n\rightarrow l}$ (by Claim 2) and there exists a $r_{1j}\in
PR_n$ such that $x_k=x_lr_{1j}.$ Now we show that $j=n.$ otherwise,
we have $j\neq n\Rightarrow r_{1j}(n)=n\Rightarrow
k=x_k(n)=x_lr_{1j}(n)=x_l(n)=l,$ which contradicts $k\neq l.$ So
$i=x_k(n-1)=x_lr_{1n}(n-1)=x_l(2)\Rightarrow
n-1=x_l^{-1}(i)=2\Rightarrow n=3,$ hence Claim 3 holds.

By Claim 3, if $X\neq B_{(i)},$ then $n=3,$ which contradicts $n\geq
5.$ Hence $X=B_{(i)},$ the assertion holds.
\end{proof}

{\bf Remark.} For $n=3,4,$ one can easily check that the result of
Theorem \ref{dl-3.1} is not true. For example, in $P_3,~
N(\{(1\,2),(1\,3\,2)\})=B^{(1)}$ and
$|\{(1\,2),(1\,3\,2)\}|=|B_{(1)}|=2,$ however,
$\{(1\,2),(1\,3\,2)\}\neq B_{(1)};$ In $P_4,~
N(\{(1\,2),(1\,2)(3\,4),(1\,3\,2),(1\,3\,4\,2),(1\,4\,2),(1\,4\,3\,2)\})=B^{(1)}$
and
$|\{(1\,2),(1\,2)(3\,4),(1\,3\,2),(1\,3\,4\,2),(1\,4\,2),(1\,4\,3\,2)\}|=|B_{(1)}|=6,$
however,
$\{(1\,2),(1\,2)(3\,4),(1\,3\,2),\\
(1\,3\,4\,2),(1\,4\,2),(1\,4\,3\,2)\}\neq
B_{(1)}.$

\begin{theorem}\label{dl-3.2}
If $n\geq 5,$ then $\Aut(P_n)=L(S_n),$ where $L(S_n)$ is the left
regular representation.
\end{theorem}

\begin{proof} For $n=5,$ a Nauty \cite {Mckay} computation
shows that $|\Aut(P_5)|=120.$ Since $|\Aut(P_5)|\geq |L(S_5)|=120,$
we have $\Aut(P_5)=L(S_5).$ We proceed by induction on $n.$ Clearly
any automorphism of $P_n$ must permute the efficient dominating sets
of $P_n.$ Let ${\mathcal B}=\{B^{(i)}: i=1,2,\cdots,n\}.$ By Theorem
\ref{dl-2.6}, $\Aut(P_n)$ naturally acts on ${\mathcal B}.$ Next we
shall show that the action of $\Aut(P_n)$ on ${\mathcal B}$ is
faithful. Assume that $\phi\in \Aut(P_n)$ such that
$\phi(B^{(i)})=B^{(i)}$ for each $i\in\{1,2,\cdots,n\}.$ By Lemma
\ref{yl-2.3}, $N(B_{(i)})=B^{(i)},$ so we have
$N(\phi(B_{(i)}))=\phi(B^{(i)})=B^{(i)},\,
|\phi(B_{(i)})|=|B_{(i)}|=|B^{(i)}|.$ By Theorem \ref{dl-3.1},
$\phi(B_{(i)})=B_{(i)}$ for each $i\in\{1,2,\cdots,n\}.$ Hence
$\phi$ can be treated as an automorphism of $P_n[B_{(n)}]=P_{n-1},$
that is, the restriction
$\phi{\upharpoonright}{B_{(n)}}\in\Aut(P_{n-1})=L(S_{n-1})$ by
induction. For the identity element $I\in S_n,$ set $y=\phi(I),$
then $y,I\in B_{(n)}^{(1)}\subseteq B_{(n)}$ and
$\phi{\upharpoonright}{B_{(n)}}=L(y).$ Hence
\begin{eqnarray*}
\phi(I)=y&\Rightarrow&\phi(N(I)\cap B_{(n)}^{(i)})=N(y)\cap
B_{(n)}^{(i)}\\
&\Rightarrow&L(y)(r_{1,i})=yr_{1,y^{-1}(i)}\\
&\Rightarrow&yr_{1,i}=yr_{1,y^{-1}(i)}\\
&\Rightarrow&y(i)=i,
\end{eqnarray*}
where $i=2,3,\cdots,n.$ So we have $\phi(I)=y=I,$ that is, $\phi$
fixes $I.$ Since $\phi(B^{(i)}_{(j)})=B^{(i)}_{(j)}$ for each
$i,j\in\{1,2,\cdots,n\},$ by Lemma \ref{yl-2.3} (ii) and the
connectedness of $P_n,$ $\phi$ fixes all vertice of $P_n,$ so
$\phi=1,$ which implies that the action of $\Aut(P_n)$ on ${\mathcal
B}$ is faithful. Thus $Aut(P_n)\lesssim \Sym({\mathcal
B})\Rightarrow |Aut(P_n)|\leq n!.$ On the other hand,
$|Aut(P_n)|\geq |L(S_n)|=n!.$ Hence $Aut(P_n)=L(S_n).$ The assertion
holds.\end{proof}

{\bf Remark.} If $n=3,$ then $P_3=C_6,$ so $\Aut(P_3)=D_{12},$ where
$D_{12}$ is the dihedral group of order $12.$ If $n=4,$ a Nauty
computation shows that $|\Aut(P_4)|=48,$ so $L(S_4)$ is a normal
subgroup of $\Aut(P_4).$ By Godsil \cite{Godsil}, $\Aut(P_4)$ is the
semiproduct $L(S_4)\rtimes\Aut(S_4,PR_4),$ where
$\Aut(S_4,PR_4)=\{\phi\in
\Aut(S_4):\,\phi(PR_4)=PR_4\}=\{1,c((2\,3))\},$ here we denote by
$1$ the identity automorphism and by $c((2\,3))$ the automorphism
induced by the conjugacy of $(2\,3)$ on $S_4.$

\end{document}